\documentclass[a4paper,10pt]{report}

\usepackage[latin1]{inputenc}
\usepackage{latexsym,amsfonts,amssymb,amsmath,amscd,euscript}

\newtheorem{thm}{Theorem}[section]

\newtheorem{dfn}[thm]{Definition}

\newcommand\sudda[1]{}

\DeclareMathOperator{\Hom}{Hom}

\newcommand\free{{\mathcal F}}
\newcommand\M{{\mathcal M}}

\newcommand\lieg{{\mathfrak g}}

\newcommand\ideal[1]{\langle#1\rangle}
\newcommand\gr{$\mathbb Z^+\times\mathbb Z_2$-graded }

\DeclareMathOperator{\ho}{H}

\DeclareMathOperator{\pp}{P}

\DeclareMathOperator{\tor}{Tor}
\DeclareMathOperator{\ext}{Ext}
\DeclareMathOperator{\der}{Der}
\DeclareMathOperator{\im}{im}

\DeclareMathOperator{\nass}{NA}
\DeclareMathOperator{\lie}{Lie}
\DeclareMathOperator{\env}{U}
\DeclareMathOperator{\ass}{Ass}
\DeclareMathOperator{\sym}{S}
\DeclareMathOperator{\ab}{Ab}
\DeclareMathOperator{\logg}{LOG}
\DeclareMathOperator{\grass}{gr}

\date{}
\title{Hilbert series, Poincar\'{e} series and homotopy Lie algebras of graded algebras }
\author{Clas L\"ofwall}
\begin{document}
\maketitle
\chapter*{Acknowledgement}
I am grateful to Rasoul Ahangari Maleki who invited me to give two lectures in the Commutative Algebra Webinars at IPM, Tehran, on February 11th and 18th, 2021. The content of the two lectures
is given in the two chapters of this note.

\chapter{Hilbert series and Poincar\'{e} series of graded algebras }

\begin{abstract}We begin with proving a formula relating the Hilbert series of a graded algebra $A$ and the Poincar\'{e} series
for $A$ in two variables using the existence of a graded minimal free resolution of the field $k$ over $A$. This gives 
the Fr\"oberg formula in the case where the bigraded $\tor^A(k,k)$ is concentrated on the diagonal, which we take as definition
of $A$ being ``Koszul". We look at a resolution in the commutative case obtained from the Koszul complex in the ``trivially Golod" case, 
which means that there are cycles that represent the homo\-logy of the Koszul complex and multiply to zero.  
The algebra structure of $\ext_A(k,k)$ is introduced in different ways. Its subalgebra generated by the one-dimensional elements is 
by definition the ``Koszul" dual of $A$.
This is studied by means of the dual of the bar-complex. We define the ``generalized Koszul complex" and construct a minimal resolution in the case 
where the cube of the augmentation ideal of $A$ is zero. This is  similar to the  trivially Golod case above, but $A$ does not need to be commutative.
The above results are at least 45 years old and most of it can be found in my thesis \cite{lof}. 
\end{abstract}

\section{Basic definitions and series}
\subsection{Basic definitions}
We will study rings $R$, that are graded algebras over a field $k$ of characteristic zero, $R=\oplus_{i=0}^\infty R_i$, where $R_iR_j\subseteq R_{i+j}$ for all $i$ and $j$. 
We also assume $R$ is connected, 
i.e., $R_0=k$,  and $R_i$ is a finite dimensional vector space over $k$ for each $i$ and that $R$ is generated by $R_1$ as an algebra. 
This means that $R$ can be presented as $R=k\ideal{x_1,\ldots,x_n}/I$ where $x_1,\ldots,x_n$ is a basis for $R_1$ and $I$ is a two-sided ideal in the free associative algebra $k\ideal{x_1,\ldots,x_n}$ (which is the same as the tensor algebra $T(R_1)$). In general $R$ may be non-commutative. It is only in section \ref{golsec} that we assume that $R$ is commutative. In this case, $R$ can be presented as a quotient of the polynomial ring $k[x_1,\ldots,x_n]$. 

We will also study the exterior algebra $\wedge(x_1,\ldots,x_n)$, which is defined as $k\ideal{x_1,\ldots,x_n}$ modulo the two-sided ideal generated by $x_ix_j+x_jx_i$ for all $i$ and $j$. 
This algebra is an example of a graded algebra with an additional $\mathbb Z_2$-grading. For an $\mathbb N\times \mathbb Z_2$-graded algebra $R$ we will call the first grading the ``degree" and the second grading the ``sign" and the sign of an element $r\in R$ will be denoted by $|r|$, which is 0 or 1, and will create a sign by $(-1)^{|r|}$. We will also use ``even" and ``odd" for elements of sign 0 and 1 respectively (the variables in the exterior algebra are odd). This kind of grading is essential when differential graded algebras are considered. 

The tensor product of two  $\mathbb N$-graded vector spaces $A$ and $B$ is defined as 
$A\otimes B=\oplus_n(A\otimes B)_n$ where $(A\otimes B)_n=\oplus_{i+j=n}A_i\otimes B_j $, and similarly for bi-graded or multi-graded vector spaces.
In case $A$ and $B$ are $\mathbb N\times\mathbb Z_2$-graded algebras, so is $A\otimes B$. 
The product of two elements $a_1\otimes b_1$ and $a_2\otimes b_2$ is
defined as $(-1)^{|b_1||a_2|}a_1a_2\otimes b_1b_2$ (and the product is extended by linearity). The exterior algebra may be seen as the tensor product of polynomial algebras on one odd variable and the polynomial algebra is the tensor product of polynomial algebras on one even variable. The free associative algebra is the ``free" product of polynomial algebras on one (even or odd) variable. The free product of two graded algebras $A$ and $B$, $A*B$, is defined as follows. Let $A_+$ and $B_+$ be the elements of positive degree in $A$ and $B$, then $(A*B)_+$ is the direct sum of all $A_+\otimes B_+\otimes A_+\ldots$ and $B_+\otimes A_+\otimes B_+\ldots$. The product is given by tensor product where two consecutive elements are multiplied if they belong to the same algebra. 

\vspace{6pt}
\noindent{\bf Example}
\begin{align*}
\text{ If }A&=k\ideal{x_1,\ldots,x_n}/J_A \text{ and }B=k\ideal{y_1,\ldots,y_m}/J_B\text{ then}\\
A*B&=k\ideal{x_1,\ldots,x_n,y_1,\ldots,y_m}/J_A+J_B
\end{align*}

For the exterior algebra, the sign of an element is the degree modulo 2. However, this does not need to be true in general, an example is $k[x_1,\ldots,x_n]\otimes\wedge(y_1,\ldots,y_n)$, where $x_1y_1$ is an odd element of degree 2.
\subsection{Series}
To any locally finite dimensional $\mathbb N$-graded vector space $V$ we associate the formal power series $V(z)=\sum_{n=0}^\infty|V_n|z^n$  
and the series for the tensor product of two such vector spaces is
$$
(U\otimes V)(z)=U(z)V(z)
$$
When $R$ is a graded algebra, we will use the notation $\ho_R(z)$ instead of $R(z)$ and we will call this series the Hilbert series of $R$. The series for the vector space 
$\tor^R(k,k)$, graded by the homological degree,  is called the Poincar\'{e} series of $R$ and denoted by $\pp_R(z)$. Since $R$ is graded $\tor^R(k,k)$ is bi-graded
(the first degree is the homological degree), and we get also a Poincar\'{e} series in two variables
$$
\pp_R(x,y)=\sum_{n,j=0}^\infty|\tor^R_{n,j}(k,k)|x^ny^j
$$
Here are some examples. In the first six examples, the Fr\"oberg formula (\ref{frob}) holds. For the last three examples we have at least the formula (\ref{form}) and the Poincar\'{e} series for these examples are obtained by the result in section \ref{r3sec}. In the first four examples a minimal resolution of $k$ over $R$ is also given. The Hilbert series for the examples may be obtained by simple combinatorial reasoning or by using the product formula for tensor products in the first two examples and the following formula for the series of a free product in the third example.
$$
1/(A*B)(z)=1/A(z)+1/B(z)-1
$$
\noindent{\bf Examples} 
\begin{align}
&R=k[x_1,\ldots,x_n]\\
&\ho_R(z)=1/(1-z)^n\notag\\
&\pp_R(z)=(1+z)^n\notag\\
&Y=k[x_1,\ldots,x_n]\otimes\wedge(T_1,\ldots,T_n),\ dT_i=x_i,\ d\text{ is extended as a derivation}\notag\\
\notag
\\
&R=\wedge(x_1,\ldots,x_n)\\
&\ho_R(z)=(1+z)^n\notag\\
&\pp_R(z)=1/(1-z)^n\notag\\
&Y=\wedge(x_1,\ldots,x_n)\otimes k[T_1,\ldots,T_n],\ dT_i=x_i ,\ d\text{ is extended as a derivation}\notag\\
\notag
\\
&R=k\ideal{x_1,\ldots,x_n}\\
&\ho_R(z)=1/(1-nz)\notag\\
&\pp_R(z)=1+nz\notag\\
&Y=k\ideal{x_1,\ldots,x_n}\otimes k[T_1,\ldots,T_n]/(T_1,\ldots,T_n)^2,\ dT_i=x_i\notag\\
\notag
\\
&R=k[x_1,\ldots,x_n]/(x_1,\ldots,x_n)^2\\
&\ho_R(z)=1+nz\notag\\
&\pp_R(z)=1/(1-nz)\notag\\
&Y=k[x_1,\ldots,x_n]/(x_1,\ldots,x_n)^2\otimes k\ideal{T_1,\ldots,T_n}\notag\\
&d(T_i|y_1|\ldots|y_k)=x_i(y_1|\ldots|y_k)\notag\\
\notag
\\
&R=k\ideal{x_1,x_2,x_3}/(x_1x_2,x_2x_3)\label{free}\\
&\ho_R(z)=1/(1-3z+2z^2-z^3)\notag\\
&\pp_R(z)=1+3z+2z^2+z^3\notag\\
\notag
\\
&R=k[x_1,x_2,x_3]/(x_1x_2,x_2x_3)\label{golex}\\
&\ho_R(z)=(1-2z^2+z^3)/(1-z)^3\notag\\
&\pp_R(z)=(1+z)^3/(1-2z^2-z^3)\notag\\
\notag
\\
&R=k[x_1,x_2,x_3]/(x_1x_2,x_2x_3)+(x_1,x_2,x_3)^3\label{r3comm}\\
&\ho_R(y)=1+3y+4y^2\notag\\
&\pp_R(x,y)=(1+xy)^2/(1-xy-x^2y^2-5x^2y^3-4x^3y^4)\notag\\
\notag
\\
&R=k\ideal{x_1,x_2,x_3}/(x_1x_2,x_2x_3)+(x_1,x_2,x_3)^3\label{r3noncomm}\\
&\ho_R(y)=1+3y+7y^2\notag\\
&\pp_R(x,y)=(1+3xy+2x^2y^2+x^3y^3/(1-16x^2y^3-11x^3y^4-7x^4y^5)\notag\\
\notag
\\
&R=k\ideal{x_1,x_2}/(x_1x_2,x_1^2,x_2^3,x_2^2x_1)\label{r3}\\
&\ho_R(y)=1+2y+2y^2\notag\\
&\pp_R(x,y)=(1+xy)/(1-xy-2x^2y^3)\notag
\end{align}

\section{A formula}
We will use the existence of a minimal graded free resolution of $k$ over $R$, where $R$ may be non-commutative. In homological degree $n$, this must be of the form $R\otimes\tor_n^R(k,k)$ with a differential that maps $\tor_n^R(k,k)$ to $R_+\otimes\tor_{n-1}^R(k,k)$. Since $R$ is connected, it can be proved by induction that $\tor_{n,j}^R(k,k)=0$ if $j<n$. Hence, in internal degree $j>0$ we have the following exact sequence.
$$
0\to\tor_{j,j}^R(k,k)\to(R\otimes\tor_{j-1}^R(k,k))_j\to\ldots\to(R\otimes\tor_{1}^R(k,k))_j\to R_j\to0
$$
This gives for each $j>0$
$$
\sum_{n=0}^j(-1)^n|(R\otimes\tor_{n}^R(k,k))_j|=0
$$
and since $|(R\otimes\tor_{0}^R(k,k))_0|=1$ we get
$$
\sum_{n,j=0}^\infty(-1)^ny^j\sum_{p+q=j}|R_p||\tor_{n,q}^R(k,k)|=1
$$
and hence
$$\sum_{j=0}^\infty\sum_{p+q=j}(|R_p|\sum_{n=0}^q(-1)^n|\tor^R_{n.q}(k,k)|)y^j$$
which gives the formula
\begin{align}\label{form}
\ho_R(y)\pp_R(-1,y)=1
\end{align}
The formula means that the coefficient for $y^n$ in the series $1/\ho_R(y)$ is the same as the alternating sum of the coefficients of the polynomial in 
the variable $x$, which is the 
coefficient for $y^n$ in $\ho_R(x,y)$. In example (\ref{r3}) the Poincar\'{e} series in two variables has the following beginning
$$
1 + 2 x y + 
 2 x^2 y^2 + (2 x^2 + 2 x^3) y^3 + (6 x^3 + 2 x^4) y^4 + (10 x^4 + 
    2 x^5) y^5 + (4 x^4 + 14 x^5 + 2 x^6) y^6 
$$
while $1/\ho_R(y)$ looks as follows
$$
1 - 2 y + 2 y^2 - 4 y^4 + 8 y^5 - 8 y^6
$$
\begin{dfn} The ring $R$ is called {\em Koszul} if $\tor_{n,j}^R(k,k)=0$ for $j\neq n$.
\end{dfn}
If $R$ is Koszul, then we have $\pp_R(x,y)=\pp_R(xy)$, and hence formula (\ref{form}) gives the Fr\"oberg formula
\begin{align}\label{frob}
\ho_R(z)\pp_R(-z)=1
\end{align}
so in this case $\pp_R(z)$ is completely determined by $\ho_R(z)$.

\section{The Golod case}\label{golsec}
In this section we assume that $R$ is commutative. 
The Koszul complex for $R$ is $K_R=R\otimes\wedge(T_1,\ldots,T_n)$, where the homological degree of $T_i$ is 1 for each $i$.
The differential is defined as the derivation which extends $dT_i=x_i$ for each $i$, where $x_1,\ldots,x_n$ is a basis for $R_1$.
This is the first step to a free minimal algebra resolution of $k$ over $R$. If the homology of $K_R$ has all Massey products defined,
then the process may be interrupted and a free minimal resolution can be obtained by tensoring $K_R$ with a free algebra. For simplicity,
we will assume that the Massey products are trivially defined, in the sense that there are cycles in $(K_R)_+$ which represent 
a basis for the homology such that 
all products of these cycles are zero.

\vspace{6pt}
\noindent{\bf Example.} Let $R=k[x_1,x_2,x_3]/(x_1x_2,x_2x_3)$ as in example (\ref{golex}) above. Then $x_1T_2,x_3T_2,x_2T_1T_3$ are cycles, which can
be proved to represent a basis for the homology of $K_R$ of positive degree and they multiply to zero.

\vspace{6pt}
Consider now the graded vector space $C=s\ho_+(K_R)$, where $s$ adds 1 to the homological degree. Hence $C_{r+1,j}=\ho_{r,j}(K_R)$ for $r>0$ and all $j$. 
Using the representing cycles for $\ho_+(K_R)$ we may define a map $\alpha:C\to Z_+(K_R)$ of homological degree -1,
 by sending a basis element to the corresponding representing cycle.
(in the Tate resolution, this is done only for 
$\ho_1(K_R)$ but here all of $\ho_+(K_R)$ is considered). A complex $Y=(K_R\otimes T(C),d)$ is defined as follows (the sign degree is the homological degree modulo 2),
where $x\in K_R$ and $c_i\in C$ for $i=1,\ldots,k$. We will use bar-notation to avoid confusion when $Y\otimes C$ is studied below.
\begin{align*}
d(x|c_1| c_2|\ldots| c_k)&=dx| c_1| c_2|\ldots| c_k +\\
&(-1)^{|x|}x\alpha(c_1)| c_2| c_3|\ldots| c_k
\end{align*} 
The square of $d$ is zero since 
\begin{itemize}
\item $d^2x=0$ for $x\in K_R$
\item $(-1)^{|dx|}(dx)\alpha(c_1)+(-1)^{|x|}(dx)\alpha(c_1)=0$ for $x\in K_R$
\item $\alpha(c_1)\alpha(c_2)=0$
\end{itemize}
We will now prove that $Y$ is acyclic (the proof is due to Gerson Levin) and hence a minimal free resolution of $k$ over $R$. 

The complex $K_R$ is a sub complex of $Y$ and as a space, 
$$Y=K_R\oplus (K_R\otimes T(C)\otimes C)=K_R\oplus (Y\otimes C)$$. 
Consider now $Y\otimes C$
as a complex with differential $d_Y\otimes 1$ and define a map $f: Y\to Y\otimes C$
 by $f(x)=0$ if $x\in K_R$ and $f(y|c)=y\otimes c$ if $y\in Y$. Then $f$ commutes with the differentials, since $d_Y(y|c)=(d_Y(y)|c)+x$, where $x\in K_R$.

Hence, we have an exact  sequence of complexes 
$$
0\to K_R\to Y\xrightarrow{f}Y\otimes C\to0
$$
which gives a long exact homology sequence
$$
\ldots\to\ho_{r+1}(Y\otimes C)\to\ho_r(K_R)\to\ho_r(Y)\to\ho_r(Y\otimes C)\to\ldots
$$
Since $\ho_0(K_R)=k$ also $\ho_0(Y)=k$. Suppose $r\geq1$ and $\ho_i(Y)=0$ for all $1\leq i\leq r-1$. 

We claim that $\ho_r(Y)=0$. Since $C_i=0$ for $i<2$, $\ho_r(Y\otimes C)=(\ho(Y)\otimes C)_r=C_r$. But $d_Y$ restricted to $C_r$ 
is $\alpha$, which is a 
monomorphism. Hence there are no non-zero cycles that are mapped to $C_r$ under the map $\ho_r(f):\ho_r(Y)\to\ho_r(Y\otimes C)$ and hence this map 
is zero. The map
$\ho_r(K_R)\to\ho_r(Y)$ is also zero, since $Y$ contains $C_{r+1}$ which kills the image of $\ho_r(K_R)$.
It follows that $\ho_r(Y)=0$, which proves the claim.

\vspace{6pt}
From the minimal free resolution $Y$, we may deduce the following formula for the Poincar\'{e} series ($n=\dim_k(R_1)$)
\begin{align}\label{gol}
\pp_R(z)=(1+z)^n/(1-z(\ho(K_R)(z)-1))
\end{align}
It is also possible to give the Poincar\'{e} series in two variables.
$$\pp_R(x,y)=(1+xy)^n/(1-x(\ho(K_R)(x,y)-1)),\ \ n=\dim(R_1)$$
In the example above we have $\ho(K_R)(z)=1+2z+z^2$ and hence we get the formula for $\pp_R(z)$ given in (\ref{golex}).
We also have
$$\ho(K_R)(x,y)=1+2xy^2+x^2y^3$$
and hence
$$\pp_R(x,y)=(1+xy)^3/(1-2x^2y^2-x^3y^3)$$
from which it follows that $R$ is Koszul, which we also know to be true by Fr\"oberg's result that any polynomial ring modulo quadratic monomials is Koszul, \cite{fro}.

\vspace{6pt}
There is a formula, similar to (\ref{form}) relating $\ho_R(y)$ and $\ho(K_R)(x,y)$. Indeed, for each internal degree $j\geq0$, the Koszul complex is the following (for $j=0$ the complex is just
$0\to k\to 0$), where
$\wedge_q$ denotes the piece of $\wedge(T_1,\ldots,T_n)$ of homological and internal degree $q$. 
$$
0\to\wedge_j\to R_1\otimes\wedge_{j-1}\to\ldots\to R_{j-1}\otimes\wedge_1\to R_j\to0
$$
We now use the fact that the Euler characteristic of a complex is the same as the Euler characteristic of the sequence of homology groups. 
This gives for each $j\geq0$
$$
\sum_{q=0}^j(-1)^{j-q}|R_q||\wedge_{j-q}|=\sum_{q=0}^j(-1)^q|\ho(K_R)_{q,j}|
$$
and hence, by multiplying with $z^j$ and summing from $j=0,\ldots,n$,
\begin{align}\label{kos}
\ho_R(z)(1-z)^n=\ho(K_R)(-1,z)
\end{align}
We check the formula in the example (\ref{golex}) studied above. We have $\ho_R(z)=(1-2z^2+z^3)/(1-z)^3$ and since $x_1T_2,x_3T_2,x_1T_2T_3$
are representing cycles for the Koszul homology, we have $\ho(K_R)(x,y)=1+2xy^2+x^2y^3$ in accordance with the formula.

\vspace{6pt}

\section{The Ext-algebra}
There is a graded algebra structure on $\ext_R(k,k)$, which we now will describe ($R$ does not need to be commutative). This algebra is generated by $\ext_R^1(k,k)$ if and 
only if $R$ is Koszul, which we will come back to in the next section. We will give four definitions of the product, but we will not prove 
that they are equal (or rather that the algebras are isomorphic up to a sign and the order of the multiplication).
\begin{enumerate}
\item Yoneda. 

An element in $\ext^n_R(k,k)$ may be considered as an equivalence class of exact sequences of $R$-modules
$$
0\to k\to A_1\to A_2\to\ldots\to A_{n-1}\to A_n\to k\to0
$$
The equivalence relation is generated by the relation $E\to E'$ for two sequences as above defined by a sequence of maps $\alpha_i$, $i=1,\ldots,n$
such that the following diagram is commutative
$$
\begin{matrix}0&\to\ k&\to\ 
A_1&\to\ A_2&\to\ldots\to&
A_{n-1}&\to\ A_n&\to\ k&\to\ 0 \\
\\
&{\scriptstyle id}\downarrow &{\scriptstyle \alpha_1}
\downarrow &{\scriptstyle\alpha_2}
\downarrow&&{\scriptstyle\alpha_{n-1}}\downarrow&
{\scriptstyle\alpha_n}\downarrow&{\scriptstyle id}\downarrow&\\
\\
0&\to\ k&\to\ A_1'&\to\ A_2'&\to\ldots\to
&A_{n-1}'&\to\ A_n'&\to\ k&\to\ 0
\end{matrix}
$$
The product of two sequences $E\in \ext_R^n(k,k)$ and $E'\in\ext_R^m(k,k)$  gives an element in $\ext_R^{n+m}(k,k)$ and is formed by putting the sequences together.
$$
0\to k\to A_1\to\ldots\to A_n\to A_1'\to\ldots\to A_m'\to k\to0
$$
Here the map  $A_n\to A_1'$ is defined by composing the maps $A_n\to k$ and $k\to A_1'$. The product is associative but in general not commutative. For more details we refer to \cite{mit}.
\item Composing chain maps.

Let $P_*$ be a projective resolution of $k$ over $R$. The complex $\Hom(P_*,P_*)$ is defined as follows. The degree $-n$ part, where $n\geq0$ (which will also be called upper degree $n$), is 
$$
\prod_{m\geq0}\Hom_R(P_{m+n},P_{m})
$$
and the differential $d$ is defined for $f$ of upper degree $n$ as
$$
(df)(x)=d(f(x))-(-1)^nf(dx).
$$
The cycles are chain maps (of a certain degree) and a boundary is a chain map, which is homotopic to zero.  Thus the cohomology of 
$\Hom(P_*,P_*)$ in upper degree $n$ consists of homotopy classes of chain maps $P_*\to P_*$ of upper degree $n$. Moreover, the map $P_*\to k$
induces a map $\Hom(P_*,P_*)\to \Hom(P_*,k)$, which is an isomorphism in cohomo\-logy (this follows by a spectral sequence argument applied to a
complete filtration). Hence the cohomology of $\Hom(P_*,P_*)$ is $\ext_R(k,k)$. The composition of chain maps endows the complex $\Hom(P_*,P_*)$ with a structure
of a DG algebra, so the cohomology has a graded algebra structure.

This definition of the product may be used to prove that $T(C^*)$ is a sub algebra of $\ext_R(k,k)$ in the Golod case (and the case $R_3=0$).
Here is a  ``lifting" formula $\tilde f:Y\to Y$ of $f\in C^*$ : 
$$\tilde f(x|c_1|c_2|\ldots|c_k)=f(c_k)x|c_1|\ldots|c_{k-1}$$
\item Free model.

Any graded algebra $R$ has a free model, which is a differential graded algebra $(T(V),d)$ together with a map $(T(V),d)\to R$, which induces an 
isomorphism in homology ($R$ is considered to have differential zero and concentrated in homological degree zero). It is even possible to
get a minimal free model (since $R$ is connected), which means that $dV\subseteq V^{\otimes\geq2}T(V)$.
This is obtained in the following way by a step by step procedure, similar to the construction of a minimal
free resolution of $k$ over $R$. 

Start to define $V_0=R_1$ and define $f:T(V_0)\to R$ as the natural epimorphism and define $d$ as zero on $V_0$. Let $K_0$ be the kernel of $f$ and define 
$V_1=s(K_0/R_+K_0)$
and choose a section $d: V_1\to K_0$ and define $f$ as zero on $V_1$. This gives that $f:T(V_0\oplus V_1)\to R$ induces an isomorphism in homological degree zero.  
Let $K_1$ be the cycles of homological degree one in $(T(V_0\oplus V_1),d)$ and define $V_2=s(K_1/R_+K_1)$ and choose a section $V_2\to K_1$. This gives that 
$f:T(V_0\oplus V_1\oplus V_2)\to R$ induces an isomorphism in homological degree 0 and 1. For the example (\ref{free}) the process stops at this point, but in general the process may continue ad infinitum. 

In example (\ref{free}) we may choose $\{T_1,T_2,T_3\}$ as a basis for $V_0$, $\{S_1,S_2\}$ as a basis for $V_1$ and $\{U\}$ as a basis for $V_2$ and define 
\begin{align*}
&f(T_i)=x_i, \ f(S_i)=f(U)=0,\\
&d(T_i)=0, \ d(S_1)=T_1\otimes T_2, \ d(S_2)=T_2\otimes T_3 \text{ and }\\
&d(U)=S_1\otimes T_3-T_1\otimes S_2.  
\end{align*} 
 \vspace{3pt}
Now $\ext_R(k,k)=sV^*$ where $V^*$ is the vector space dual of $V$ and $s$ adds 1 to the cohomological degree and changes sign.
 The product is defined as the dual of the quadratic part $d_2$ of $d$, i.e., 
$V^*\otimes V^*\xrightarrow{d_2^*}V^*$. 
For the example (\ref{free}) we get $\ext_R^{i+1}(k,k)=V_i^*$ for $i=0,1,2$ and $\ext_R^i(k,k)=0$ for $i\geq4$ and hence the global homological dimension of $R$ is 3 and the Poincar\'{e} series is $1+3z+2z^2+z^3$ as stated in (\ref{free}). The  product is given by 
$$T_1^*T_2^*=S_1^*,\ T_2^*T_3^*=S_2^*, \ S_1^*T_3^*=U^*,\ T_1^*S_2^*=U^*$$
 and all other products are zero. The sign in the last product follows the rule given in the next section.

It is easy to prove e.g., that the multiplication is associative.
\item Cobar construction

The bar resolution $B$ is of the form $(R\otimes T(R_+),d)$ and is a free resolution of $k$ over $R$. Hence $\ext_R(k,k)$ is the cohomology of $\Hom(B,k)$, which is
$(T(R_+^*),d^*)$, where the cohomological degree is the length of a tensor ($R_+^*$ has cohomological degree 1). The induced differential on $k\otimes B=T(R_+)$ is the alternating sum of terms where two consecutive elements in a tensor product are multiplied. From this it is easy to see that $d^*$ is a derivation on $T(R_+^*)$ which extends the dual of the multiplication map $R_+\otimes R_+\to R_+$. It follows that $(T(R_+^*),d^*)$ is a differential graded algebra and hence its cohomology is an algebra. This definition of the product on $\ext_R(k,k)$ will be used in the next section. Observe that $(T(R_+^*),d^*)$ is not a minimal model of $\ext_R(k,k)$ in general. There is no map in general from the model to $\ext_R(k,k)$, which induces an isomorphism in homology. If there is such a map, then the double $\ext$-algebra would be $R_+$, which implies that $R$ is Koszul, and indeed if $R$ is Koszul, then the map exists, making $(T(R_+^*),d^*)$ a minimal model of $\ext_R(k,k)$.
\end{enumerate}

\section{The Koszul dual}
\begin{dfn} The Koszul dual of $R$ is defined as $R^!=\oplus_{n\geq0}\ext_R^{n,n}(k,k)$, which is a subalgebra of $\ext_R(k,k)$.
\end{dfn}
 We will prove that $R^!$ is generated by 
$\ext_R^{1,1}(k,k)$ and has quadratic relations. The diagonal elements, i.e., the elements where the cohomological degree is equal to the internal degree, in $T(R_+^*)$ are of the form $v_1\otimes v_2\otimes\ldots\otimes v_n$ where $v_i\in R_1^*$ for all $i$. Since the image of the multiplication map $R_+\otimes R_+\to R_+$ is in $R_{\geq2}$,  it follows that these elements are cycles and since $R_1^*$ generates $T(R_1^*)$ under tensor product we get that $R^!$ is generated by $\ext_R^{1,1}(k,k)$. The only possibility for an element in $T(R_1^*)$ to be a boundary is that it is a boundary of an element in $R_{i_1}^*\otimes R_{i_2}^*\ldots\otimes R_{i_n}^*$, where there is a $j$ such that $i_r=1$ for all $r\neq j$ and $i_j=2$. If 
$$
\varphi: R_1\otimes R_1\to R_2
$$ 
is multiplication, then the boundaries in $T(R_1^*)$ of cohomological degree $n$ are 
$$
\sum_{i=0}^{n-2}(R_1^*)^{\otimes i}\otimes \im(\varphi^*)\otimes (R_1^*)^{\otimes(n-i-2)}
$$
which is precisely the two-sided ideal in $T(R_1^*)$ generated by $\im(\varphi^*)$. Thus $R^!$ has quadratic relations
$$
R^!=T(R_1^*)/(\im(\varphi^*))
$$
 and these relations may be computed from the quadratic relations in $R$, which is $\ker(\varphi)$. We have $\im(\varphi^*)=(\ker(\varphi))^0=$ the elements of $R_1^*\otimes R_1^*$ which are zero on $\ker(\varphi)$. These elements may be found by solving a system of homogeneous  linear equations. 

\vspace{6pt}
\noindent{\bf Example}
$$
R=k\ideal{x_1,x_2}/(x_1x_2+x_2x_1+x_2^2,x_1x_2+x_1^2)
$$
Let $\{T_1,T_2\}$ be the dual basis of $\{x_1,x_2\}$ and consider $S=aT_1T_2+bT_2T_1+cT_1^2+dT_2^2$. Then $S(x_1x_2+x_2x_1+x_2^2)=a+b+d$ and 
$S(x_1^2+x_1x_2)=a+c$. A basis for the solutions of the system 
\begin{align*}
a+b+d&=0\\
a+c&=0
\end{align*}
is $(1,-1,-1,0),(0,-1,0,1)$. Hence
$$
R^!=k\ideal{T_1,T_2}/(T_1T_2-T_2T_1-T_1^2,T_2^2-T_2T_1)
$$

Since the (co)homological degree of $R$ is zero and
the cohomological degree of $R_+^*$ is one in $(T(R_+^*),d^*)$, there should be a sign shift going from $R$ to $R^!$. If the variables in $R$ are even, 
then the variables in $R^!$ should be odd and vice versa. We have assumed in the example above that $(f\otimes g)(u\otimes v)=f(u)g(v)$. 
This works if the variables in $R$ have equal signs, but if there are both even and odd variables, then the following rule should be used
$$
(f\otimes g)(u\otimes v)=(-1)^{|v|(|u|+1)}f(u)g(v),
$$
see \cite{frolof}.

 \vspace{6pt}
\noindent{\bf Example}
Consider $R=k[x]\otimes\wedge(y)$, where $x$ is even and $y$ is odd. Then $R^!$ should be $\wedge(x^*)\otimes k[y^*]$, where $x^*$ is odd and $y^*$ is even.
We have $R=k\ideal{x,y}/(xy-yx,y^2)$ and if we forget the signs, the relations in $R^!$ would be $x^*y^*+y^*x^*,(x^*)^2$, but we should have $x^*y^*=y^*x^*$. 
Using the new rule instead, we get $(x^*y^*)(xy)=-1$ 
and $(y^*x^*)(yx)=1$ and hence $(x^*y^*-y^*x^*)(xy-yx)=0$, which gives the correct relations in $R^!$.

\vspace{6pt}
Taking the Koszul dual of $R^!$ gives $(R^!)^!=T(R_1)/(\ker(\varphi))$, which is easily seen using the multiplication map on $R^!$ 
$$
\phi: R_1^*\otimes R_1^*\to (R_1^*\otimes R_1^*)/\im(\varphi^*)
$$
and using $\im(\phi^*)=(\im(\varphi^*))^0=\ker(\varphi.)$
\section{The generalized Koszul complex}
The vector space dual of $R^!$ in  degree $n$ may be seen to be
$$
K_n=\bigcap_{i=0}^{n-2}(R_1)^{\otimes i}\otimes \ker(\varphi)\otimes (R_1)^{\otimes(n-i-2)}\text{ for }n\geq2
$$
and $K_0=k$, $K_1=R_1$.
Here it is used that $(V/U)^*=U^0$, $(U_1+U_2)^0=U_1^0\cap U_2^0$, $V^{**}=V$ and $(\im(\varphi^*))^0=\ker(\varphi)$.
If $B$ is the bar resolution, then $K_n$ is the  subspace of $B\otimes k$ consisting of the diagonal elements such that all terms of the differential are zero. 
In this way, 
$R\otimes K_n$ defines a free minimal sub complex of $B$, the {\em generalized Koszul complex}, 
$$GK_R=R\otimes (R^!)^*$$
Here is the complex in low homological degrees: 
$$R\otimes (R_1\otimes \ker(\varphi)\cap\ker(\varphi)\otimes R_1)\to R\otimes\ker(\varphi)\to R\otimes R_1\to R$$
We have
\begin{align}\label{exsec}
0\to K_{n+1}\to R_1\otimes K_n\to R_2\otimes K_{n-1}
\end{align}
 is exact for $n\geq0$.
 Hence, cycles in $GK_R$ are in $R_{\geq2}\otimes(R^!)^*$ plus boundaries.
 
$GK_R$ is acyclic if and only if $R$ is Koszul. Even if $R$ is not Koszul one may in some occasions use $GK_R$ to build a minimal resolution of $k$ over $R$, see next section and \cite{roos}.

Let $Y=R\otimes\tor^R(k,k)$ be a minimal graded free resolution of $k$ over $R$. By degree reason, the differential $d_Y$ maps $\tor^R_{n,n}(k,k)$ to
$R_1\otimes\tor^R_{n-1,n-1}(k,k)$ and thus $R\otimes(\oplus_n\tor^R_{n,n}(k,k))$ is a sub complex of $Y$, which may be identified with $GK_R$.

We also have
$$
GK_{R^!}=(GK_R)^*
$$
 which can be useful in computing the homology of $GK_R$, see \cite{roos}, Appendix B, e.g., Theorem B.9.
\section{The case $R_3=0$}\label{r3sec}
We will now use $GK_R$ to obtain a minimal free resolution of $k$ over $R$ in a situation similar to the Golod case in section \ref{golsec}. If $R$ is commutative,
then the method is precisely the same as in section \ref{golsec}, since in this case $GK_R$ is a differential graded algebra. 
In the non-commutative case we have to assume that the representing cycles for  $\ho(GK_R)$ multiplied with boundaries are zero. This is fulfilled if $R_3=0$, which we will assume to be true in the sequel.

Let $C=s\ho_+(GK_R)$. Hence $C_0=C_1=0$ and $C_{n,j}=\ho_{n-1,j}(GK_R)$ for $n\geq2$ and all $j$. Using (\ref{exsec}) we have for each $n\geq0$ the following exact sequence
\begin{align}\label{r3exsec}
0\to K_{n+1}\to R_1\otimes K_n\to R_2\otimes K_{n-1}\to C_n\to 0
\end{align}
It follows that $C_n$ is concentrated in internal degree $n+1$ for $n\geq2$. Let  $\alpha_n$ be a section $\alpha_n:C_n\to R_2\otimes K_{n-1}$ for $n\geq2$.
We now define a differential $d$ on $Y=GK_R\otimes T(C)$ by
\begin{align*}
d(x\otimes y)&=dx\otimes y\text{  for }x\in (GK_R)_+,\ y\in T(C)\\
d(x\otimes y\otimes z)&=x\alpha_n(y)\otimes z\text{ for }x\in R,\ y\in C_n,\ z\in T(C)
\end{align*}
We have $d^2=0$ since $(dx)\im(\alpha_n)=0$ for $x\in(GK_R)_1$. The proof that $Y$ is acyclic is the same as the proof given in section \ref{golsec}. 
 Let $C(z)=\sum_{n\geq0}|C_n|z^n$ and $C(x,y)=\sum_{x,y\geq0}x^ny^j|C_{n,j}|$. Since $C_{n,j}=0$ for $j\ne n+1$, we get $C(x,y)=yC(xy)$ and since 
 $R^!$ is concentrated
 on the diagonal, we have $\ho_{R^!}(x,y)=\ho_{R^!}(xy)$. Hereby we get the Poincar\'{e} series for $R$ as
 $$
 P_R(x,y)=\ho_{R^!}(xy)/(1-yC(xy))
 $$
 We may use (\ref{r3exsec}) to compute $C(z)$. The alternating  sum of dimensions of the sequence (\ref{r3exsec}) is zero. Multiply with $z^{n+1}$ and take the sum for $n\geq0$ to get
 $$
 \sum_{n\geq0}z^{n+1}|R^!_{n+1}|-|R_1|z\sum_{n\geq0}z^n|R^!_n|+|R_2|z^2\sum_{n\geq1}z^{n-1}|R^!_{n-1}|-zC(z)=0
 $$
 Hence
 $$
 zC(z)=\ho_{R^!}(z)-1-|R_1|z\ho_{R^!}(z)+|R_2|z^2\ho_{R^!}(z)
 $$
 which finally gives
 $$
 P_R(x,y)=\frac{x\ho_{R^!}(xy)}{1+x-\ho_{R^!}(xy)(1-|R_1|xy+|R_2|x^2y^2)}
 $$
 The series $\ho_{R^!}(z)$ may be computed up to a certain degree using the package {\em Bergman} developed by J\"orgen Backelin, see \cite{ber}.
 
 \vspace{10pt}
 The elements in the $\ext$-algebra off the diagonal is always a two-sided ideal in  $\ext_R(k,k)$ and the quotient algebra is $R^!$. Hence,
 $R^!$ is always a direct factor of $\ext_R(k,k)$ as algebras.


\chapter{The homotopy Lie algebra of a graded commutative algebra}
\begin{abstract}
We give the definition of a graded Lie algebra with examples. The free Lie algebra on a set $X$ is defined and
the enveloping algebra of quotients of free Lie algebras are studied. The Koszul dual is looked upon as the enveloping algebra of a Lie algebra in the graded commutative case with examples. 
The Poincar\'{e}-Birkhoff-Witt theorem is stated and as a consequence, a formula for the Hilbert series of the enveloping algebra is obtained. In the graded commutative case a coproduct $\Delta$
on $\ext_A(k,k)$ is possible to  define and thereby the homotopy Lie algebra is defined as the set of $x$ such that
$\Delta(x)=x\otimes1+1\otimes x$. 
Other definitions of the homotopy Lie algebra are given
by means of minimal algebra resolutions and models.  The Lie subalgebra $\eta$
generated by the one-dimensional elements in the homotopy Lie algebra $\lieg$ of $A$ is studied.  Here, the enveloping algebra of $\eta$ is the Koszul dual of $A$.
Examples of homotopy Lie algebras are given for
complete intersections, Golod rings and rings with the cube of the augmentation ideal equal to zero.  The dimensions of a Lie algebra given the Hilbert series for its enveloping algebra can be computed by means of a logarithmic formula (see \cite{lof4}, section 2). We  apply this to edge ideals. Some examples of periodic Lie subalgebras $\eta$ are given, yielding examples of irrational Poincar\'{e} series.  The ``holonomy" Lie algebra (see \cite{lof5},\cite{lof3}) of a hyperplane arrangement is defined, and some examples are given, such as the graphical arrangement $K_4$. 
\end{abstract}
\section{Definition and examples}
A positively graded Lie super algebra (called just a Lie algebra in the sequel) is a \gr vector space $L$ with a bilinear operator $[\ ,\ ]$ satisfying the axioms ``anti-symmetry" and ``Jacobi identity" given below. 
The $\mathbb Z_2$-grading is called the ``sign" and 
the sign of an element $a$ is denoted by $|a|$, which is 0 or 1, and will create a sign by $(-1)^{|a|}$. We will also use ``even" and ``odd" for elements of sign 0 and 1 respectively. 

The $\mathbb Z^+$-degree will be called just the degree and the sign of an element are often equal to the degree modulo 2, but not always. 
\begin{align*}
[a,b]&=-(-1)^{|a||b|}[b,a]\quad&\text{anti-symmetry}\\
[a,[b,c]]&=[[a,b],c]+(-1)^{|a||b|}[b,[a,c]]\quad&\text{Jacobi identity}
\end{align*}
There are extra axioms in characteristic 2 and 3, e.g., $[x,[x,x]]=0$ for odd $x$ in characteristic 3, but we will assume as always that the characteristic of the field is zero.
A Lie algebra $L$ is called  ``ordinary" if $|a|=0$ for all $a\in L$.
Observe that Jacobi identity may also be written
$$
[[a,b],c]=[a,[b,c]]+(-1)^{|b||c|}[[a,c],b].
$$
The main example of a Lie algebra is a \gr associative algebra with Lie product $[a,b]=ab-(-1)^{|a||b|}ba$. In fact, as we will see, any Lie algebra is a Lie sub algebra of this kind of Lie algebra. 

A derivation $D$ on a \gr associative or non-associative algebra $A$ is a graded $k$-linear map $A\to A$ satisfying
$$
D(ab)=D(a)b+(-1)^{|D||a|}aD(b)
$$
Jacobi identity expresses the fact that $[a,\cdot]$ is a graded derivation on the Lie algebra.
It is easy to prove that if $D_1,D_2$ are derivations, then  $[D_1,D_2]=D_1D_2-(-1)^{|D_1||D_2|}D_2D_1$ is also a derivation. Thus the set of derivations on $A$, $\der(A,A)$ is a
Lie algebra, but the degree may be any integer. To get a $\mathbb Z^+$-graded Lie algebra, one may restrict to only positive degrees or only negative degrees, the latter is used 
when $\der_R(Y,Y)$ is studied for a minimal free $R$-algebra resolution $Y$, which we will come back to. Here $\der_R(Y,Y)$ is the set of derivations $D$ on $Y$ such that $D(r)=0$ for $r\in R$.

Another important example is $\free(X)$, the free Lie algebra on a \gr set $X$. It is defined by first consider the set $B$ of all binary non-associative ``monomials" in the alphabet $X$, where the binary operation is written as $[\ ,\ ]$. An element in $B$ is e.g., $[[x,y],[x,[x,y]]]$ where $x,y\in X$. Each monomial has a degree and a sign. Now define 
$$
\nass(X)=\oplus_{n,s}\nass(X)_{n,s},
$$
where $\nass(X)_{n,s}$ is the vector space with basis consisting of all monomials in $B$ of degree $n$ and sign $s$. The operation $[\ ,\ ]$ is extended by bilinearity to an operation on 
homogeneous elements of $\nass(X)$. In this way $\nass(X)$ is a non-associative algebra. Let $I$ be the two-sided ideal in $\nass(X)$ generated by all elements of the form
$[a,b]+(-1)^{|a||b|}[b,a]$ and $[a,[b,c]]-[[a,b],c]-(-1)^{|a||b|}[b,[a,c]]$ for homogeneous elements $a,b,c\in\nass(X)$. Finally $\free  (X)$ is defined as
$$
\free  (X)=\nass(X)/I.
$$
Obviously, $\free  (X)$ is a Lie algebra and also any Lie algebra which is generated by $X$ is a quotient of $\free  (X)$, which justifies the use of the name ``free Lie algebra" for $\free  (X)$. 

Classical matrix Lie algebras obtained from Lie groups, such as $sl(n,\mathbb C)$ are not graded in our sense, but they will still come into play 
as the period in a periodic Lie algebra. 

\section{The enveloping algebra}
We have seen that any \gr associative algebra $A$ defines a Lie algebra. This defines a functor which we call $\lie(A)$. There is a left adjoint to $\lie$ called $\env$, the enveloping algebra functor. This means that for
any Lie algebra $L$, there is a Lie algebra map $\epsilon_L:L\to\lie(\env(L))$ such that for any Lie algebra $L$ and associative algebra $A$ and a Lie algebra map $f:L\to \lie(A)$ 
there is a unique algebra map $g:\env(L)\to A$ such that $\lie(g)\circ\epsilon_L=f$.

Any Lie algebra has an enveloping algebra. The free Lie algebra $\free  (X)$ has the enveloping algebra $\ass(X)$, which is the tensor algebra on the vector space spanned by $X$. Let $\epsilon_X$ be the map from $\free  (X)$ to $\lie(\ass(X))$ induced by the identity map on $X$. If $I$ is a Lie ideal in $\free  (X)$ then the enveloping algebra of $\free  (X)/I$ is $\ass(X)/(\epsilon_X(I))$. All this is easily proved.

\noindent{\bf Example}
\begin{align*}
L&=\free  (a,b,c)/([a,b]+[c,c],[a,c]+[b,c])\quad a,b\text{ even and }c\text{ odd}\\
\env(L)&=k\ideal{a,b,c}/(ab-ba+2c^2,ac-ca+bc-cb)
\end{align*}

In the previous lecture we learned that the Koszul dual of a quadratic algebra $R=k\ideal{x_1,\ldots,x_n}/I$ is $R^!=k\ideal{T_1,\ldots,T_n}/I^0$, where $T_1,\ldots,T_n$ is the dual basis
of $x_1,\ldots,x_n$ and $I^0$ is generated by the set of quadratic elements in $k\ideal{T_1,\ldots,T_n}$ which are zero on the generators of $I$. If $R$ is commutative (in the ordinary sense), then $I$ contains the commutators $x_ix_j-x_jx_i$ for all $i$ and $j$. This means that the relations in $R^!$ are linear combinations of $T_iT_j+T_jT_i$ for all $i,j$ and hence $R^!$ is the enveloping algebra of a Lie algebra, $L_R$, which is a quotient of $\free  (T_1,\ldots,T_n)$, where $T_1,\ldots,T_n$ are odd, modulo the Lie ideal generated by the quadratic relations in $R^!$ 
(considered as elements in $\free  (T_1,\ldots,T_n)$).

\vspace{6pt}
 \noindent{\bf Example}
\begin{align*}
R=\ &k[x_1,x_2,x_3]/(x_1x_2,x_2x_3)\\
R^!=\ &k\ideal{T_1,T_2,T_3}/(T_1^2,T_2^2,T_3^2,T_1T_3+T_3T_1)=\\
&k\ideal{T_1,T_2,T_3}/([T_1,T_1],[T_2,T_2],[T_3,T_3],[T_1,T_3])\\
L_R=\ &\free  (T_1,T_2,T_3)/([T_1,T_1],[T_2,T_2],[T_3,T_3],[T_1,T_3])\quad T_i \text{ odd for }i=1,2,3
\end{align*}

\vspace{6pt}
If $R$ is skew commutative, then the ideal of $R$ contains all $x_ix_j+x_jx_i$ and hence $R^!$ has relations which are linear combinations of $T_iT_j-T_jT_i$ for $i<j$, which means that $R^!$ is the enveloping algebra of an ordinary Lie algebra (the signs are zero).

\vspace{6pt}
 \noindent{\bf Example}
\begin{align*}
R=\ &\wedge(x_1,x_2,x_3)/(x_1x_2+x_2x_3)\\
R^!=\ & k\ideal{T_1,T_2,T_3}/(T_1T_2-T_2T_1-T_2T_3+T_3T_2,T_1T_3-T_3T_1)=\\
&k\ideal{T_1,T_2,T_3}/([T_1,T_2]-[T_2,T_3],[T_1,T_3])\\
L_R=\ &\free  (T_1,T_2,T_3)/([T_1,T_2]-[T_2,T_3],[T_1,T_3])\quad T_i \text{ even for }i=1,2,3
\end{align*}
One has to be careful in computing the Koszul dual when there are both even and odd variables.
We have 
$$(x^*y^*)(xy)=-1 \text{ if } |x|=0,\ |y|=1,\text{ otherwise }(x^*y^*)(xy)=1$$
The Koszul dual of 
$$R=k\ideal{x,y}/(xy-yx,y^2)=k[x]\otimes\wedge(y), \ |x|=0,\ |y|=1$$ is 
$$R^!=k\ideal{x^*,y^*}/(x^*y^*-y^*x^*,(x^*)^2)=\wedge(x^*)\otimes k[y^*], \ |x^*|=1,\ |y^*|=0$$
 and $R^!$ is the enveloping algebra of 
 $$\free(a,b)/([a,b],[a,a]), \ |a|=1,\ |b|=0
 $$

\vspace{6pt}
Here is the example from above

\vspace{6pt}
 \noindent{\bf Example}
\begin{align*}
L&=\free  (a,b,c)/([a,b]+[c,c],[a,c]+[b,c])\quad a,b\text{ even and }c\text{ odd}
\end{align*}
We have that $L=L_R$ for the following algebra $R$
 \begin{align*}
R&=\wedge(x,y)\otimes k[z]/(2xy-z^2,xz-yz)\quad x,y\text{ odd and }z\text{ even}
\end{align*}

In the examples above $R$ is Koszul, so $R^!=\ext_R(k,k)$ and the Lie algebras $L_R$ given in the examples are in fact the homotopy Lie algebra of $R$. In general when
$R$ is commutative or skew commutative or a mixture of these cases, the Lie 
algebra $L_R$ given by $R^!=\env(L_R)$ is the Lie sub algebra of the homotopy Lie algebra generated by the elements of degree one.  

\section{The Poincar\'{e} Birkhoff Witt theorem}
The PBW theorem states that a $k$-basis for $\env(L)$ is obtained from an ordered $k$-basis $B$ of $L$ as all ordered monomials of basis elements from $\epsilon_L(B)$, except that
the monomials do not contain repetition of odd basis elements (this is because a square of an odd element $x$ belongs to the Lie algebra, as $\frac{1}{2}[x,x]$). In particular it follows 
that the map $\epsilon_L:L\to \env(L)$ is injective (we drop from now on the Lie functor), and from this it follows that $\free  (X)$ is isomorphic to the Lie sub algebra of $\ass(X)$ generated by $X$.  

\vspace{6pt}
\noindent{\bf Example}
 \begin{align*}
L&=\free  (a,b)/([a,a],[b,b])\quad a,b\text{ odd}
\end{align*}
An ordered k-basis for $L$ is $a<b<[a,b]$. 

A $k$-basis for $\env(L)$ is $ababa\ldots,babab\ldots$. This is however not a PBW-basis. The ordered monomials with no 
repetitions of $a$ and $b$ are  $a^\alpha b^\beta[a,b]^i$, where $i\geq0$ and $\alpha,\beta=0,1$.

It follows from the PBW theorem that the series for $U(L)$ is the same as the series for the symmetric algebra on $L$. In fact, there is a filtration on $\env(L)$ such that the 
associated graded ring is isomorphic to $\sym(L)$. Thus we get the product formula given below for the series of $\env(L)$ in terms of the series for $L$. 


In the general case $L(z)$ is the sum of two series, one for the even elements $L_{even}(z)$ and one for the odd elements $L_{odd}(z)$. If
the sign of an element
in $L$ is the same as the degree modulo 2, then $L_{even}$ is the series with even $z$-exponents, and $L_{odd}$ is the series with odd $z$-exponents.

\vspace{6pt}
Suppose
\begin{align*}
&L_{even}(z)=\sum_{i=1}^\infty\alpha_iz^i\\
&L_{odd}(z)=\sum_{i=1}^\infty\beta_iz^i
\end{align*}

then
\begin{align}\label{prod}
\env(L)(z)=\prod_{i=1}^\infty\frac{(1+z^{i})^{\beta_{i}}}{(1-z^{i})^{\alpha_{i}}}
\end{align}

\section{The homotopy Lie algebra}
We have seen that if $R$ is (skew) commutative and Koszul, then $\ext_R(k,k)$ is the enveloping algebra of a Lie algebra, which makes it possible to define the homotopy Lie
algebra of $R$ as this Lie algebra. It is true in general, in the (skew) commutative case, that $\ext_R(k,k)$ is the enveloping algebra of a Lie algebra. This follows from a theorem by Milnor-Moore, Andr\'{e} and Sj\"odin. The key point is that there should exist a comultiplication  on $\ext_R(k,k)$, $\Delta$,  which is an algebra homomorphism $\Delta:\ext_R(k,k)\to\ext_R(k,k)\otimes\ext_R(k,k)$ ($\Delta$ is the dual of a multiplication on $\tor^R(k,k)$)). The Lie algebra is then obtained as the set of ``primitive" elements, i.e., elements $x\in\ext_R(k,k)$ satisfying $\Delta x=x\otimes1+1\otimes x$. 

There are different ways to ``construct" the homotopy Lie algebra without knowing in beforehand $\ext_R(k,k)$. We present here three different methods.
\begin{enumerate}
\item Cobar construction

The tensor algebra $T(R_+^*)$ has a comultiplication $\Delta$ which is defined such that it is multiplicative with respect to the tensor product and such that the elements in
$R_+^*$ are primitive. (You may check that if $a,b$ are primitive then also $[a,b]=ab-(-1)^{|a||b|}$ is primitive. We have seen that $T(R_+^*)$ is the enveloping algebra of
the free Lie algebra $\free  (R_+^*)$.) One can prove that $\Delta$ commutes with the differential $d^*$ and thus defines a comultiplication on the cohomology which is $\ext_R(k,k)$. The method to find elements in the homotopy Lie algebra should therefore be to find primitive cycles (modulo boundaries) in $(T(R_+^*),d^*)$ and Lie multiply them by means of the 
rule $ab-(-1)^{|a||b|}ba$.
\item Commutative Models

A (skew) commutative ring $R$ has a minimal model, which is a polynomial algebra on even and odd generators of homological degrees $\geq0$ together with a map to $R$ which induces an isomorphism in homology ($R$ is assumed to be concentrated in homological degree zero). The construction is similar to the free minmal model which was studied in the previous lecture.

\vspace{6pt}
\noindent{\bf Example}
\begin{align*}
&f:k[T_1,T_2,S_1,S_2; dT_i=0,dS_1=T_1^2,dS_2=T_2^3]\to 
k[x_1,x_2]/(x_1^2,x_2^3)\\
&\text{where }\deg(T_i)=0\text{ and }T_i \text{ is even, }\\
&\deg(S_i)=1\text{ and }S_i \text{ is odd and }\\
&f(T_i)=x_i,\ f(S_i)=0
\end{align*}

\vspace{6pt}
The homotopy Lie algebra is obtained as the dual of the space of variables in the minimal model, where the homological degrees of the variables are raised by one and also
the signs are changed. The Lie product is obtained as the dual of the quadratic part of the differential (there is no linear part since the model is minimal), more precisely,  
\begin{align*}
\text{If \ } dU&=\sum_{S<T}\lambda_{STU}ST+\sum_{|S|=0}\lambda_{SU}S^2+ \text{ higher terms} \\
\text{then \ }[S^*,T^*]&=-\sum_U(-1)^{|S^*||T|}\lambda_{STU}U^* \text{ \ if }S<T\text{ and}\\
(S^*)^2&=\frac{1}{2}[S^*,S^*]=-\sum_U\lambda_{SU}U^*\text{ \ if }|S|=0
\end{align*}

In the example above we get $(T_1^*)^2=-S_1^*$ and all other products are zero.
\item Minimal algebra resolutions

If $Y$ is a free minimal algebra resolution of $k$, then $\der_R(Y,Y)$ is a Lie algebra as we saw above. It is also a sub complex of $\Hom(Y,Y)$, since
the differential $\delta$ on $Y$ is a derivation and hence if $D\in \der_R(Y,Y)$ then also $[\delta,D]\in \der_R(Y,Y)$ 
and $d(D)=[\delta,D]$ is the differential of $D$ as an element in $\Hom(Y,Y)$. 
Moreover, $\der_R(Y,Y)$  is a differential Lie algebra in the sense that $d$ is a derivation on $\der_R(Y,Y)$, which follows by Jacobi identity (observe that $d=[\delta,\cdot]$).
Now we may define the homotopy Lie algebra, $L_R$, of $R$, as 
$$
L_R=\ho(\der_R(Y,Y)).
$$
Observe that we only consider ``negative" degrees (which we rather want to call positive upper degrees), in the sense that a derivation $D$ of upper degree $n\geq1$ maps $Y_{m+n}$
to $Y_m$ for all $m\geq0$.

 For each variable $S$ in $Y$, there is a derivation $j_S$ which maps $S$ to 1 and all variables of 
degree less than $S$ to zero,  and it commutes with the differential (i.e., it is a cycle in $\der_R(Y,Y)$). 
These maps considered as cycles in $CC=\Hom_R(Y,Y)$ generate, by composition, $\ext_R(k,k)$ as an algebra (Sj\"odin, \cite{sj1}) 
and they form a $k$-basis for the Lie algebra.
Thus as a vector space $L_R$ may be identified with the space of variables in $Y$,
but contrary to the case with minimal models above, there is no change of homological degree and sign. The Lie product is given by $[j_S,j_T]=j_S\circ j_T-(-1)^{|S||T|}j_T\circ j_S$.

The same example as above looks as follows.

\vspace{6pt}
\noindent{\bf Example}
\begin{align*}
Y=&R[T_1,T_2,S_1,S_2; dT_i=x_i,dS_1=x_1T_1,dS_2=x_2^2T_2]\\
&\text{where } \deg(T_1)=1,\ T_i \text{ odd, }\deg(S_i)=2,\ S_i \text{ even}\\
&j_{T_1}(S_1)=-T_1,\ j_{T_1}(S_2)=0,\ j_{T_2}(S_1)=0,\ j_{T_2}(S_2)=-x_2T_2.
\end{align*}
It follows that $(j_{T_1})^2=-j_{S_1}$, all other products are zero ($(j_{T_2})^2$ is a boundary, since it is mapped to zero by the augmentation $Y\to k$).
 
\end{enumerate}
\section{The homotopy Lie algebra of some classes of rings}
\noindent{\bf complete intersections}

We have seen above an example of the Lie algebra of a complete intersection. In general, for a complete intersection $R$, the homotopy Lie algebra $L_R$ is zero in degrees $\geq3$. It is generated by $(L_R)_1$ if and only if
the relations in $R$ are quadratic. Each minimal relation of degree $>2$ in $R$ corresponds to a minimal generator for $L_R$ of degree 2. The quadratic relations in $L_R$ are obtained by Koszul duality from the quadratic relations in $R$ and the quadratic algebra defined by the quadratic relations in $R$ is Koszul. Let $\eta_R$ be the Lie sub algebra of $L_R$ generated by $(L_R)_1$. Then $L_R$ may be described as the direct sum of $\eta_R$ and an abelian Lie algebra on generators of degree 2 corresponding to the minimal relations of degree $>2$ in $R$.

\vspace{6pt}
\noindent{\bf Golod rings}

If $R$ is a Golod ring, then the Lie sub algebra of $L_R$ consisting of all elements of degree $\geq2$ is free on the generator set $s\ho_+(K_R)$. This is an ideal in $L_R$ and the quotient is the abelian Lie algebra on $(L_R)_1$, $\ab((L_R)_1)$. Hence, there is an exact sequence of Lie algebras (which is non-split),
$$
0\to\free  (s\ho_+(K_R))\to L_R\to \ab((L_R)_1)\to0
$$
One may also consider the Lie sub algebra of $L_R$, $J=\oplus_{q>p}(L_R)_{pq}$, of all elements off the diagonal. It is a Lie sub algebra of $(L_R)_{\geq2}$ and hence free, since it is a general fact that a Lie sub algebra of a free Lie algebra is free (Shirshov-Witt-Shtern theorem). We have that $J$ is an ideal of $L_R$ and the quotient is $\eta_R$. Hence, there is a split (to the right) exact sequence 
$$
0\to J\to L_R\to\eta_R\to0
$$  

\vspace{6pt}
\noindent{\bf Rings with $R_3=0$ }

This case is similar to the Golod case. Also here the Lie sub algebra of $L_R$, $J=\oplus_{q>p}(L_R)_{pq}$, is free. The generators of $J$ as a Lie algebra are the elements on the super  diagonal, $C^*=\oplus_{p\geq2}(L_R)_{p,p+1}$ and $C^*$ may be obtained as the third syzygy of a minimal free resolution of $k$ over $R^!$:
$$
0\to C^*\to R_2^*\otimes R^!\to R_1^*\otimes R^!\to R^!\to k\to0
$$
Hence, we have the following exact sequence of Lie algebras, which is split to the right
$$
0\to \free  (C^*)\to L_R\to\eta_R\to0
$$
\section{A logarithmic formula}
Given the series for an enveloping algebra of a Lie algebra $L$, it is possible, using the product formula (\ref{prod}), to compute recursively the dimensions of $L$ in different degrees. We will give a closed formula for $L(z)$ in terms of the series $\env(L)(z)$, which generalizes Serre's formula for the dimensions of a free Lie algebra on $d$ even generators.
The formula is slightly more difficult when the signs are general, so we first consider the cases when either all signs are zero or the sign of an element is the degree modulo 2.
Suppose $V(z)$ is a series with rational coefficients and constant term 1, then $\logg_0(V)(z)$ and $\logg_1(V)(z)$ are series with rational coefficients and constant term 0, defined by
$$
\logg_0(V)(z)=\sum_{r=1}^\infty\frac{\mu(r)}{r}\log(V(z^r))
$$
and
$$
\logg_1(V)(z)=\sum_{r=1}^\infty\frac{\mu(r)}{r}\log(V((-1)^{r+1}z^r))
$$
where $\mu$ is the M\"obius $\mu$-function
It follows that $\logg_i(V_1V_2)=\logg_i(V_1)+\logg_i(V_1)$ for $i=0,1$. We now claim that if $V(z)$ is the series for the enveloping algebra of a Lie algebra $L$ concentrated in sign zero, and $\dim(L_i)=\epsilon_i$, i.e., if
 $$
V(z)=\prod_{i=1}^\infty\frac{1}{(1-z^{i})^{\epsilon_{i}}},
$$
then 
$$
\logg_0(V)(z)=\sum_{i=1}^\infty\epsilon_iz^i
$$
We have that $\logg_0$ commutes with the substitution $z\to z^i$ for any $i\geq1$. Hence, it is enough to prove the claim for $\logg_0$ in the only case when $V(z)=1/(1-z)$. But
$$
\sum_{r=1}^\infty\frac{\mu(r)}{r}\log(1/(1-z^r))=\sum_{r,s=1}^\infty\frac{\mu(r)}{rs}z^{rs}=\sum_{n=1}\frac{1}{n}z^n\sum_{r|n}\mu(r)=z
$$
since $\sum_{r|n}\mu(r)=0$ for $n>1$ and $\mu(1)=1$. 

For the second case, we claim that if $V(z)$ is the series for the enveloping algebra of a Lie algebra $L$, where the sign is equal to the degree modulo 2, and $\dim(L_i)=\epsilon_i$, i.e., if
$$
V(z)=\prod_{i=1}^\infty\frac{(1+z^{2i+1})^{\epsilon_{2i+1}}}{(1-z^{2i})^{\epsilon_{2i}}}
$$ 
then
$$
\logg_1(V)(z)=\sum_{i=1}^\infty\epsilon_iz^i
$$
It is in fact enough to prove the claim for $V(z)=1+z$ and $V(z)=1/(1-z^2)$ and this is done in the same way as above.

\vspace{6pt}
Here is the general case  when the sign of the Lie algebra is not necessarily the degree modulo 2.
Suppose $V_0(z)$ and $V_1(z)$ are the even and odd series for a \gr vector space $V$.
Consider $V(z,y)=V_0(z)+yV_1(z)$ as an element in the ring $\mathbb Z[[z]][y]/(y^2-1)$.
Then define
$$\logg(V)(z,y)=\sum_{r=1}^\infty\frac{\mu(r)}{r}\log(V(z^r,(-1)^{r+1}y^r))$$
If $$V(z,y)=\prod_{i=1}^\infty\frac{(1+yz^{i})^{\beta_{i}}}{(1-z^{i})^{\alpha_{i}}}$$ then
$$\logg(V)(z,y)=\sum_{i=1}^\infty\alpha_iz^i+y\sum_{i=1}^\infty\beta_iz^i$$ 

\vspace{6pt}
The free Lie algebra $L$ on $d$ even generators has $\env(L)(z)=1/(1-dz)$. Applying $\logg_0$, one gets Serre's formula:
$$
L(z)=-\sum_{r=1}^\infty\frac{\mu(r)}{r}\log(1-dz^r)=\sum_{r,s=1}^\infty\frac{\mu(r)}{rs}d^sz^{rs}=\sum_{n=1}^\infty z^n\frac{1}{n}\sum_{r|n}\mu(n/r)d^r
$$

\vspace{6pt}
An application is for Koszul algebras. Suppose $R$ is Koszul. Then by Fr\"oberg's formula,
$$L_R(z)=\logg_1(1/\ho_R(-z))=\sum_{r,s=1}^\infty\frac{\mu(r)}{rs}(1-\ho_R((-z)^r))^s$$
In the previous lecture, we proved the following formula expressing the Hilbert series in terms of the homology of the Koszul complex ($n=\dim(R_1)$).
$$
\ho_R(z)(1-z)^n=\ho(K_R)(-1,z)
$$
Let $d_R(z)=\ho(K_R)(-1,z)$. If $R$ is Koszul we get by Fr\"oberg's formula
$$
\pp_R(z)=\frac{(1+z)^n}{d_R(-z)}
$$
 Using the logarithmic formula $\logg_1$ above, we get the following formula for $L_R(z)$.
 $$
 L_R(z)=nz+\sum_{r,s=1}^\infty\frac{\mu(r)}{rs}(1-d_R((-z)^r))^s
 $$
 In particular, this formula holds for a ring $R$ defined by an edge ideal.
 \section{Periodic Lie algebras}
 The homotopy Lie algebra of a ring $R$ cannot be periodic, since by  a theorem of F\'{e}lix-Halperin-Thomas, the Lie algebra $L_R$ grows faster than any polynomial if $R$ is not a complete intersection. However, the  Lie sub algebra generated by $(L_R)_1$, $\eta_R$, may be periodic. In fact, my student Anna Larsson proved that the ``periodization" of any finite dimensional simple contragredient (Kac's terminology) Lie super algebra $\mathfrak g$ except $sl(2)$ has quadratic relations (and generators of degree one). The periodization is $\mathfrak g\otimes tk[t]$, where $t^i$ has degree $i$ and $[xt^i,yt^j]=[x,y]t^{i+j}$. David Anick found the odd generated Lie algebra $\free(a,b)/([a,[a,[a,b]]],[b,[b,[b,a]]])$ with series $2,3,2,1,2,3,2,1,\ldots$. It is a periodization of $sl(2|1)$ (Kac notation). To obtain an example of a Lie algebra which comes from a commutative algebra, the relations have to be quadratic. The first such example was found by Jan-Erik Roos and myself \cite{Lo-Ro}.  We started with a quadratic algebra $R$ in 5 variables, such that $\eta_R$ seemed to be periodic by experiment. At last we could prove that it was formed as a periodization of a Lie algebra which occurs in Kac's list of simple finite dimensional Lie super algebras. This gave a new example of an algebra $R$ with $R_3=0$ and $\pp_R(z)$ irrational. The first such example was found by David Anick.
 
 I will show you an example of a periodization of the ordinary Lie algebra $sl(3)$ equipped with a $\mathbb Z_2$-grading, with the property that the relations are quadratic.  The Lie algebra $sl(3)$ consists of all $3\times3$-matrices with trace zero. The Lie product $[a,b]$ is defined as $ab-ba$. Here is the subdivision of $sl(3)=L_1\oplus L_2$. We use the notation $e_{jk}$ for the standard basis elements in a $3\times 3$-matrix.
 \begin{align*}
 L_1 \text{ has basis } &h_1=e_{11}-e_{22},h_2=e_{22}-e_{33},\\
      &a_1=e_{12}+e_{21},a_2=e_{23}+e_{32},a_3=e_{13}+e_{31}\\
 L_2\text{ has basis }&e_{12}-e_{21},e_{23}-e_{32},e_{13}-e_{31}
 \end{align*} 
Hence $L_1$ is 5-dimensional and $L_2$ is 3-dimensional. It is easy to prove that $[L_1,L_1]=L_2$, $[L_1,L_2]=L_1$ and $[L_2,L_2]=L_2$. Thus $L_1\oplus L_2$ is a $\mathbb Z_2$-grading of $sl(3)$. We now define the periodization $L$ as $\oplus_{i\geq0}L_1t^{2i+1}\oplus_{i\geq1}L_2t^{2i}$. It follows from above that $L$ is generated by $L_1t$ and it is easy to prove that the following seven quadratic relations hold among the generators of degree 1, $L_1t$ in $L$. We write $x$ instead of $xt$ for the generators.
\begin{align*}
&[h_1,a_1]+2[h_2,a_1]=0,\ [h_2,a_2]+2[h_1,a_2]=0,\ [h_1,a_3]-[h_2,a_3]=0,\\
&[h_1,a_3]-[a_1,a_2]=0,\ [h_1,a_2]+[a_1,a_3]=0,\ [h_2,a_1]-[a_2,a_3]=0,[h_1,h_2]=0
\end{align*}
The question is now whether there are more relations, or if the free Lie algebra on $L_1$ modulo these relations is periodic with the dimensions $5,3,5,3,5,3,\ldots$. 
This can be examined by a program, which is a package in Macaulay2, called GradedLieAlgebras, see \cite{lun}, \cite{mac2}.

The ring $R$ associated to this example (in the sense that $\eta_R=L$) is an exterior algebra on 5 generators $x_1,x_2,y_1,y_2,y_3$ modulo the following three quadratic relations, 
$$
2x_1y_1+y_1x_2+y_2y_3,\ \ 2x_2y_2+y_2x_1+y_1y_3,\ \ x_1y_3+x_2y_3+y_1y_2
$$
A surprising fact about this example is that the dimension sequence \hfill\break
$5,3,5,3,5,3,\ldots$ for $\eta_R$ is generic in the variety of all quadratic rings in 5 odd variables modulo 3 quadratic relations. One could think that the periodic behaviour is an odd property among all examples, but in this case this is the ``normal" case. There is another interesting phenomenon in this example, which seems to occur in every periodic example. If the field has characteristic $p$ the dimension of the Lie algebra is one higher at one place in the period and if the degree is divisible by $p$. In this case, if the field has characteristic 5, the dimension of $L$ in degrees divisible by 10 is 4 instead of 3.
\section{The holonomy Lie algebra of a matroid}
A matroid is a set of subsets $\M$ of a finite set $E$, which satisfies the following two axioms.
 \begin{align*}
 1)&\quad A\in\M, \ B\subset A\quad \Rightarrow \quad B\in\M\\
 2)&\quad A,B\in\M, \ |A|<|B|\quad\Rightarrow \quad\exists x\in B\setminus A \text{ such that }A\cup\{x\}\in\M
 \end{align*}
 The axioms are obviously fulfilled if $E$ is a finite set of vectors in a vector space and $\M$ is the set of independent subsets of $E$. For any matroid $\M$, we will use  ``$A$ independent" to mean $A\in\M$ and ``$A$ dependent" for $A\notin\M$. A central hyperplane arrangement is a finite set of hyperplanes through the origin in $\mathbb C^n$. The set $E$ of normal vectors of the hyperplanes defines a matroid as above. The complement $X$ of the union of the hyperplanes has an interesting cohomology algebra, called the Orlik-Solomon algebra (with coefficients in a field $k$ of characteristic zero).  Its structure depends only on the matroid. 
 
 The Orlik-Solomon algebra may be defined for any matroid $\M$ over \hfill\break $E=\{e_1,\ldots,e_n\}$, 
 which is ``simple", i.e., all subsets of $E$ with at most two elements belong to $\mathcal M$. (A central hyperplane arrangement is simple, since two hyperplanes through the origin cannot be parallel.)
 The Orlik-Solomon algebra, $\mathcal O$, is defined as a quotient of the exterior algebra $\wedge(e_1,\ldots,e_n)$, modulo the ideal generated by the following polynomials $p_A$, one for each dependent subset of $E$, $A=\{a_1,\ldots,a_k\}$.
 $$
 p_A=\sum_{i=1}^k(-1)^ia_1a_2\cdots \hat a_i\cdots a_k
 $$
 The holonomy Lie algebra of $\mathcal M$ is defined as 
 $\eta_\mathcal O$, i.e., the Lie sub algebra of the homotopy Lie algebra of $\mathcal O$ generated by the elements of degree one. 
 
 The quadratic relations in $\mathcal O$ are $ab-ac+bc$ for each dependent 3-set $\{a,b,c\}$. 
 To find the relations in $\eta_\mathcal O$, suppose $x_1,\ldots,x_n$ is the dual basis of the variables $e_1,\ldots,e_n$ in $\mathcal O$.  Hence, $x_1,\ldots,x_n$ are even. For each monomial $ab$ in $\mathcal O$ one has to find all 
 dependent 3-sets $\{a,b,c\}$. If there is no dependent 3-set containing $a$ and $b$, then the corresponding dual variables commute in $\eta_\mathcal O$. Here $x_i,x_j$ ``commute" means
 $[x_i,x_j]=x_ix_j-x_jx_i=0$. The union of all dependent 3-sets containing $a$ and $b$ is called a   ``2-flat", which by definition is a maximal subset of $E$ such that all 3-subsets are dependent. Two different 2-flats have at most one element in common. 
 
 Suppose $\{a,b,c\}$ is a 2-flat.
We have $ab-ac+bc=0$ in $\mathcal O$ and no other quadratic relation contains $ab$. Hence the dual relations are $[x,y]+[x,z]= [x,z]+[y,z]=0$. From these relations it follows that $x+y+z$ commutes with $x,y,z$.
 In general, for each 2-flat $A=\{e_{i_1},\ldots,e_{i_k}\}$ one gets the following relations in $\eta_\mathcal O$. Let $y=x_{i_1}+\cdots+x_{i_k}$, then 
  \begin{align}\label{rel}
 [x_{i_j},y] =0\text{ for } j=1,\ldots,k
 \end{align}
 Hence $\eta_\mathcal O$ is the free Lie algebra $\free(x_1,\ldots,x_n)$ modulo the relations (\ref{rel}) for each 2-flat $A$. 
 
 The relations in (\ref{rel}) are not linearly independent, the sum of them is zero.
 
 For each 2-flat $A=\{e_{i_1},\ldots,e_{i_k}\}$ of size at least three one may consider the ``local" Lie sub algebra, $\mathcal L_A$,  of $\eta_\mathcal O$ generated by $x_{i_1},\ldots,x_{i_k}$. We have that the element $y=x_{i_1}+\cdots+x_{i_k}$ is central in $\mathcal L_A$ and also $\mathcal L_A/y$ is free in degrees $\geq2$ on $k-1$ generators of degree one. The local Lie sub algebras 
 $\mathcal L_A$ are ``glued" together to form $\eta_\mathcal O$. Sometimes $\eta_\mathcal O$ is the direct sum of the local Lie sub algebras $\mathcal L_A$ and in this case 
 $\eta_\mathcal O$ is called ``decomposable". This happens precisely when 
 $[x,[y,z]]=0$ whenever $\{x,y,z\}$ is independent. 
 
 Suppose $\mathcal F$ is a set of subsets of a finite set $E$ satisfying that two different elements in $\free$ contain at most one element in common and all  elements in $\free$ are of size at least three. Then $\free$ defines uniquely a simple matroid of subsets of size at most three such that $\free$ is the set of all 2-flats of size at least three, see \cite{lof5}.
 
 \vspace{6pt}
 \noindent{\bf Example}
 
 Any graph defines a matroid over the set of edges by letting an independent subset be a subset with no circuits (a forest). The same matroid may also be defined as the matroid of a central hyperplane arrangement consisting of all hyperplanes $x_i-x_j$ for which there is an edge between the vertices numbered $i$ and $j$. Any 2-flat must be of size at most three, since 4 edges must contain 3 edges which do not form a triangle. In fact, the 2-flats of size three are precisely the triangles and the 2-flats of size 2 is a set of two edges with no vertex in common. Let us study $K_4$ in detail. Here is a picture with names on the 6 edges.
 $$
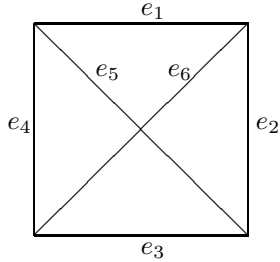
\begin{picture}(100,120)  
\put(40,103){$e_1$} 
\put(83,60){$e_2$}
\put(40,13){$e_3$}
\put(-10,60){$e_4$}
\put(23,80){$e_5$}
\put(50,80){$e_6$}
\put(0,20){\line(1,0){80}}
\put(0,20){\line(0,1){80}}
\put(80,20){\line(0,1){80}}
\put(0,100){\line(1,0){80}}
\put(0,100){\line(1,-1){80}}
\put(0,20){\line(1,1){80}}
\end{picture} 
$$
We have the following 2-flats of size three
 $$
 \{e_1,e_2,e_5\},\ \{e_3,e_4,e_5\},\ \{e_1,e_4,e_6\},\ \{e_2,e_3,e_6\}
 $$ 
 and the following 2-flats of size two
 $$
 \{e_1,e_3\},\ \{e_2,e_4\},\ \{e_5,e_6\}
 $$
 This gives the holonomy Lie algebra as the free Lie algebra
 $$
 \free(x_1,x_2,x_3,x_4,x_5,x_6),
 $$
 where the variables are even, modulo the ideal generated by
 \begin{align*}
 &[x_ 1,x_ 2+x_ 5],[x_ 2,x_ 1+x_ 5],\\
 &[x_ 3,x_ 4+x_ 5],[x_ 4,x_ 3+x_ 5],\\
 &[x_ 1,x_ 4+x_ 6],[x_ 4,x_ 1+x_ 6],\\
 &[x_ 2,x_ 3+x_ 6],[x_ 3,x_ 2+x_ 6],\\
 &[x_ 1,x_ 3],[x_ 2,x_ 4],[x_5,x_6]
 \end{align*}
 The holonomy Lie algebra for $K_4$ is not decomposable, since e.g., $[x_3,[x_4,x_6]]\neq0$. In fact, the ideal generated by $\{x_ 3,x_ 4,x_ 6\}$ is a free Lie sub algebra and the full holonomy Lie algebra is an extension of this Lie algebra and the local Lie algebra $L_{\{e_1,e_2,e_5\}}$.
 
 The holonomy Lie algebra of a central hyperplane arrangement with complement $X$ is of interest in itself. There is a result by Kohno \cite{koh}, that the holonomy Lie algebra is equal to 
 $$
 \grass(\pi_1(X))\otimes_{\mathbb Z}k
 $$
 where for a group $G$, $\grass(G)$ is the graded associated with respect to the lower central series of $G$, which is a Lie algebra over $\mathbb Z$.


\begin{thebibliography}{10}
\bibitem{av} L.L. Avramov, {Local algebra and rational homotopy}, 
Homotopie Alg\'{e}brique et Alg\`{e}bre Locale (Luminy 1982), Ast\'{e}risque 113--114, Soc. Math. France, Paris, (1984), 15--43.
\bibitem{ad} J.F. Adams, {On the non-existence of elements of Hopf invariant one}, Annals of Math, 72, (1960), 20--104, chapter 2.
\bibitem{ber} J. Backelin, {The Gr\"obner basis calculator program {\em Bergman}}, available at servus.math.su.se/bergman.
\bibitem{fro}R. Fr\"oberg, Determination of a class of Poincar\'{e} series, Math. Scand., vol 37, (1975), 29--39.
\bibitem{frolof} R.Fr\"oberg and C. L\"ofwall, Koszul homology and Lie algebras with application to generic forms and points, 
Homology, Homotopy and Applications,  vol 4(2), (2002), 227--258.
\bibitem{koh} T. Kohno, {On the holonomy Lie algebra and
the nilpotent completion of the fundamental group of the complement of
hypersurfaces}, Nagoya Math. J., 92, (1983), 21--37.
\bibitem{lof1}C. L\"ofwall, {On the homotopy Lie algebra of a local ring}, JPAA, vol 38, (1985), 305--312.
\bibitem{lof2}C. L\"ofwall, {Central elements and deformations of local rings}, JPAA, vol 91, Issues 1--3, (1994), 183--192.
\bibitem{lof5} C. L\"ofwall, {Decompositions theorems for a generalization of the holo\-nomy Lie algebra of an arrangement}, Communications in Algebra, vol 44, Issue 11, (2016), 4654--4663. 
\bibitem{lof3}C. L\"ofwall, {The holonomy Lie algebra of a matroid}, arXiv: 2012.12044.
\bibitem{lof4}C. L\"ofwall, {Cyclic homology of algebras of global dimension at most two}, arXiv: 1711.03644v2.
\bibitem{lof}C. L\"ofwall, {On the subalgebra generated by the one-dimensional elements in the Yoneda Ext-algebra}, 
in J.-E. Roos, ed.,Algebra, Algebraic Topology and their interactions, Proceedings, Stockholm 1983, 
Lecture Notes in Mathematics 1183, (1986), Springer-Verlag, Berlin, 291--338. 
\bibitem{Lo-Ro} C. L\"ofwall, J.-E. Roos,
       {A nonnilpotent $1$-$2$-presented graded Hopf algebra whose Hilbert series converges in the unit circle},Adv. Math., vol 130, (1997), 161--200.
\bibitem{lun}S.Lundqvist, C. L\"ofwall, {Software for doing computations in graded Lie algebras}, Journal of Software for Algebra and Geometry, forthcoming.
\bibitem{mac2}D. Grayson, M. Stillman, {Macaulay2, a software system for research in algebraic geometry}, Available at  https://faculty.math.illinois.edu/Macaulay2/

\bibitem{mit}B.Mitchell, {Theory of categories}, Pure and Applied Mathematics, vol 17, Elsevier, 1965.
\bibitem{roos}J.-E. Roos, {A computer-aided study of the graded Lie algebra of a local commutative noetherian ring}, 
Journal of Pure and Applied Algebra, vol 91, (1994), 255--315.
\bibitem{sj1} G. Sj\"odin, {A set of generators for $\ext_R(k,k)$}, Math. Scand, vol 38, No 2, (1976), 199--210.
\bibitem{sj2} G. Sj\"odin, {The $\ext$-algebra of a Golod ring}, JPAA, vol 38, (1985), 337--351.
\end{thebibliography}
\end{document}